\documentclass{amsproc}
\usepackage{amssymb}
\usepackage{amsbsy}
\usepackage{amscd}

\usepackage[mathscr]{eucal}
\usepackage{nicefrac}

\DeclareSymbolFont{rsfs}{U}{rsfs}{m}{n}
\DeclareSymbolFontAlphabet{\mathrsfs}{rsfs}

\usepackage[all]{xy}
\usepackage{tikz}
\usetikzlibrary{arrows,decorations.pathmorphing,backgrounds,positioning,fit,petri}
\usetikzlibrary{decorations.markings}

\usepackage{verbatim}
\usepackage{version}
\usepackage{color}
\newenvironment{NB}{
\color{red}{\bf NB}. \footnotesize
}{}

\excludeversion{NB}
\excludeversion{NB2}
\makeatletter

\usepackage{url}
\usepackage{ifmtarg}

\usepackage{xr-hyper}
\usepackage[%
unicode,
pdftex,
bookmarks=true,
colorlinks=true,
pdfnewwindow=true]{hyperref}

\usepackage{xspace}

\usepackage[all]{xy}

\usepackage{stmaryrd}
\SetSymbolFont{stmry}{bold}{U}{stmry}{m}{n}
\usepackage{pigpen} 
\usepackage{enumitem}
\usepackage{calc}

\usepackage{textgreek}

\usepackage[whole,autotilde]{bxcjkjatype}

\renewcommand{\thesubsection}{\thesection(\@roman\c@subsection)}
\newcounter{number}
\makeatother
\newtheorem{Theorem}[equation]{Theorem}

\theoremstyle{definition}

\theoremstyle{remark}

\numberwithin{equation}{section}

\newcommand{\thmref}[1]{Theorem~\ref{#1}}
\newcommand{\secref}[1]{\S\ref{#1}}
\newcommand{\subsecref}[1]{\S\ref{#1}}

\newcommand{\CC}{{\mathbb C}}
\newcommand{\ZZ}{{\mathbb Z}}

\newcommand{\RR}{{\mathbb R}}
\newcommand{\proj}{{\mathbb P}}
\newcommand{\CP}{\proj}

\newcommand{\SU}{\operatorname{\rm SU}}
\newcommand{\GL}{\operatorname{GL}}

\newcommand{\SO}{\operatorname{\rm SO}}
\newcommand{\grpSp}{\operatorname{\rm Sp}}

\newcommand{\su}{\operatorname{\mathfrak{su}}}

\newcommand{\Spec}{\operatorname{Spec}\nolimits}

\newcommand{\Hom}{\operatorname{Hom}}

\newcommand{\DB}{\overline{\partial}}

\renewcommand{\MR}[1]{}
\newcommand{\Wedge}{{\textstyle \bigwedge}}

\newcommand{\dslash}{/\!\!/}
\newcommand{\vin}[1]{\operatorname{i}(#1)} 
\newcommand{\vout}[1]{\operatorname{o}(#1)} 
\newcommand{\bM}{\mathbf M}
\newcommand{\bN}{\mathbf N}
\newcommand{\CS}{\mathrm{CS}}

\newcommand{\shfO}{\mathcal O}
\newcommand{\tslash}{/\!\!/\!\!/}
\newcommand{\tslabar}{\mathbin{
\setbox0=\hbox{/\!\!/\!\!/}\rule[0.4\ht0]{\wd0}{.3\dp0}\kern-\wd0\box0}}

\newcommand{\bmu}{\boldsymbol\mu}

\newcommand{\la}{\lambda}

\newcommand{\Gr}{\mathrm{Gr}}

\newcommand{\cR}{\mathcal R}

\newcommand{\cT}{\mathcal T}
\newcommand{\cK}{\mathcal K}
\newcommand{\cO}{\mathcal O}

\newcommand{\scP}{\mathscr P}

\makeatletter
\newcommand{\cA}[1][{}]{%
  \@ifmtarg{#1}%
  {\mathcal A}
  {\mathcal A(#1)}
}
\newcommand{\cAh}[1][{}]{%
  \@ifmtarg{#1}%
  {\mathcal A_\hbar}
  {\mathcal A_\hbar(#1)}
}
\makeatother

\newcommand{\ft}{\mathfrak t}

\newcommand{\po}{\ar@{}[dr]|{\text{\pigpenfont R}}}
\newcommand{\pb}{\ar@{}[dr]|{\text{\pigpenfont J}}}
\newcommand{\pp}{\ar@{}[dr]|{\text{\pigpenfont P}}}

\newcommand{\cM}{\mathcal M}

\begin{document}
\title[Coulomb branches of $3d$ $N=4$ gauge theories] {Introduction to a
provisional mathematical definition of Coulomb branches of $3$-dimensional
  $\mathcal N=4$ gauge theories
}
\author[H.~Nakajima]{Hiraku Nakajima}
\address{Research Institute for Mathematical Sciences,
Kyoto University, Kyoto 606-8502,
Japan}
\email{nakajima@kurims.kyoto-u.ac.jp}
\thanks{The research of the author is supported by JSPS Kakenhi Grant Numbers
24224001, 
25220701, 
16H06335. 
}
\subjclass[2000]{}

\date{}

\begin{abstract}
    This is an introduction to a provisional mathematical definition
    of Coulomb branches of $3$-dimensional $\mathcal N=4$
    supersymmetric gauge theories, studied in
    \cite{2015arXiv150303676N,main}. This is an expanded version of an
    article \cite{2016arXiv161209014N} appeared in 第61回代数学シンポジ
    ウム報告集 (2016), written originally in Japanese.
\end{abstract}

\maketitle

\setcounter{tocdepth}{2}

  \section{Coulomb and Higgs branches -- complex symplectic varieties and deformation quantization}\label{sec:CH}

Let $G$ be a complex reductive group and $\bM$ its symplectic
representation. Namely $\bM$ is a vector space with a symplectic form
$\omega$, and $G$ acts linearly on $\bM$ preserving $\omega$. Let us denote the Lie algebra of $G$ by $\mathfrak g$.

The mathematical definition of the \emph{Coulomb branch} of $3d$ SUSY
gauge theory gives a recipe to construct a complex affine-algebraic
symplectic variety\footnote{It has a singularity in general. It is
  expected that singularities is symplectic in the sense of Beauville,
  but the proof is not given.}  $\cM_C\equiv \cM_C(G,\bM)$ from
$(G,\bM)$:
\begin{equation*}
  (G,\bM) \leadsto \cM_C(G,\bM).
\end{equation*}
It is motivated by a research in a quantum field theory in physics.
It is different from known constructions of algebraic varieties, such
as zero sets of polynomials, quotient spaces, moduli spaces, etc. We
first construct the coordinate ring $\CC[\cM_C]$ as a homology group
with convolution product. Then we define $\cM_C$ as its spectrum, and
study its geometric properties.

As we will explain later, $\cM_C$ is birational to $T^* T^\vee/W$:
\begin{equation*}
    \cM_C \approx T^*T^\vee/W = \ft\times T^\vee/W.
\end{equation*}
In physics, the right hand side is regarded as the classical
description of the Coulomb branch, and $\cM_C$ is obtained from it
after \emph{quantum correction}.
Here $T^\vee$ is the dual of a
maximal torus $T$ of $G$, and $W$ is the Weyl group. $T^*T^\vee$ is
the cotangent bundle of $T^\vee$, and $\ft$ is the Lie algebra of $T$.
In particular, the birational class of $\cM_C$ depends only on $G$. It
is independent of the representation $\bM$.

As we have already mentioned above, we construct a ring as a homology
group with convolution product. This method has been used successfully
in geometric representation theory. Since study of representation is
the main motivation there, it is usual to construct a \emph{noncommutative}
algebra.
In fact, also for the Coulomb branch, we \emph{do} get a
noncommutative deformation $\cAh$ of $\cM_C$ simultaneously. Here a
noncommutative deformation means a noncommutative algebra $\cAh$
defined over $\CC[\hbar]$ such that $\cAh/\hbar\cAh$ is isomorphic to
the coordinate ring $\CC[\cM_C]$ and the Poisson bracket
\begin{equation*}
    \{ f, g\} = \left.
    \frac{\tilde f\tilde g - \tilde g\tilde f}{\hbar}\right|_{\hbar=0}, \qquad
    \tilde f|_{\hbar=0} = f, \quad \tilde g|_{\hbar=0} = g
\end{equation*}
is equal to one given by the symplectic form. We call
$\cAh\equiv\cAh(G,\bM)$ the \emph{quantized Coulomb branch}.

Many noncommutative algebras studied in representation theory are
deformation of commutative algebras, e.g., the universal enveloping
algebra $U(\mathfrak g)$ of a Lie algebra $\mathfrak g$ is a
deformation of the symmetric algebra of $\mathfrak g$. However it is
remarkable (at least to the author) that commutative algebras and its
deformation appear in a systematic construction.

In the first paper \cite{2015arXiv150303676N}, we consider a general
$\bM$, but we only constructed $\CC[\cM_C]$ as a vector space. A
definition of a product was given later in \cite{main}, under the
assumption that $\bM$ is of a form $\bM = \bN\oplus\bN^*$ (cotangent
type condition). A physical argument says that the induced
homomorphism $\pi_4(G)\to \pi_4(\grpSp(\bM))\cong\{\pm 1\}$ must
vanish in order to have a well-defined Feynman measure on the space of
fields.\footnote{This is pointed out by Witten via Braverman. It is
  possibly related to an existence of an orientation data for the
  vanishing cycle considered in \cite{2015arXiv150303676N}.} We do not
know whether this vanishing is required or enough to define a Coulomb
branch at this moment, but the assumption $\bM = \bN\oplus\bN^*$ is
too restrictive, as more general cases have been studied in physics
literature.
We will later use the notation $\cM(G,\bN)$ when we assume
$\bM = \bN\oplus\bN^*$ after \secref{sec:def}. There should be no fear
of confusion.

There is another well-known recipe to construct a complex
affine-algebraic symplectic variety from $(G,\bM)$. It is the
\emph{symplectic reduction}
\begin{equation*}
    \bM\tslash G = \mu^{-1}(0)\dslash G,
\end{equation*}
called the \emph{Higgs branch} of the same $3d$ SUSY gauge theory
associated with $(G,\bM)$ in the physics literature.
Here $\mu\colon\bM\to\mathfrak g^*$ is the moment map
vanishing at the origin, and $\mu^{-1}(0)\dslash G$ is the quotient
space of $\mu^{-1}(0)$ by $G$ in the sense of geometric invariant
theory, namely the coordinate ring $\CC[\mu^{-1}(0)\dslash G]$ is the
space of $G$-invariant polynomials $\CC[\mu^{-1}(0)]^G$ in the
coordinate ring of $\mu^{-1}(0)$.

When $\bM = \bN\oplus\bN^*$, the ring $\mathcal D(\bN)$ of polynomial
coefficients differential operators on $\bN$ gives a noncommutative
deformation of $\bM$. (In order to introduce $\hbar$, one consider the
Rees algebra associated with the filtration given by degrees of
differential operators.) A noncommutative analog of the symplectic
reduction has been known as a quantum symplectic reduction, which
should be considered as an appropriate `quotient' of $\mathcal D(\bN)$
of $G$. It gives a noncommutative deformation of $\cM_H$.

In representation theory, we have experienced that interesting
symplectic varieties and their quantization appear as symplectic
reductions, e.g., quiver varieties and toric hyper-K\"ahler manifolds.
On the other hand, the study of Coulomb branches is just started.
We get symplectic varieties, whose description as symplectic reduction
of finite dimensional symplectic vector spaces are not known. Hence we
expect the importance of Coulomb branches will increase in future.

We believe that representation theory of a quantized Coulomb branch
$\cAh$ is easier to study than that of a quantum symplectic reduction,
as it is defined as a convolution algebra, hence powerful geometric
techniques (see \cite{CG}) can be applied.

Also, the pair of Higgs and Coulomb branches of a given $(G,\bM)$ is
expected to be a symplectic dual pair in the sense of
Braden-Licata-Proudfoot-Webster \cite{2014arXiv1407.0964B} in many
cases. The symplectic duality expects a mysterious relation between a
pair of symplectic varieties. The whole picture of the symplectic
duality is yet to be explored, but it at least says that it is
meaningful and important to study Higgs and Coulomb branches
simultaneously. It should be noted that the current framework of
symplectic duality in \cite{2014arXiv1407.0964B} requires both
symplectic varieties have symplectic resolutions. This assumption is
not satisfied in many examples of Higgs and Coulomb branches. Hence we
should start to look for more general framework of the symplectic duality.

\section{Physical background}

In \secref{sec:CH} I have explained why study of Coulomb branches
could be interesting from mathematical point of view. In this section
I will try to explain physical background, as far as I can. Like
\cite{main} this article is written so that no knowledge on physics is
required to read except this section. The reader does not need to
understand this section, as I myself does not well. But my superficial
understanding led me to find a definition given in the next
section, and it is of my belief that some understanding of physics
background will be necessary to achieve new results in Coulomb
branches.
A reader in hurry could skip this section, but it is of my hope that
(s)he does not.

Let me emphasize that I, by no means, intend to ignore past research in
physics, which strongly motivated us to obtain most of results
explained in this paper. The relevant literature can be found in
\cite{2015arXiv150303676N}.

In physics like differential geometry, people use a maximal compact
subgroup $G_c$ of a complex reductive group $G$. Similarly we assume
that $\bM$ has an inner product preserved by $G_c$.

A given pair $(G,\bM)$, physicists associate a $3$-dimensional
supersymmetric gauge theory.
It is an example of quantum field theories which are defined by path
integrals of lagrangians over infinite dimensional space of all
fields.
There are two important fields, one is a connection on a principal
$G_c$ bundle $P$ over $\RR^3$, and the other is a section of $P$ with
values in $\bM$. Other fields are spinors and sections of vector
bundles associated with $P$. They play important role in physics, but
we ignore them as we will only give rough understanding. Anyhow the
lagrangian containing curvature of the connection and differential of
the section is well-defined functional, but the path integral does not
have a mathematically rigorous definition.
Configurations giving local minimum of lagrangian are classical
solutions of motion in quantum mechanics, hence are important objects.
In our situation local minimum configurations form a finite
dimensional space, instead of a single path. It is called \emph{the
  moduli space of vacua}. In fact, this will not be the right
definition, and it gives the \emph{classical} moduli space, and the actual moduli space receives \emph{corrections} as we will explain below.

The lagrangian is sum of square norm of the curvature and the
derivative of the section and others. Local minimum are attained when
several terms vanish. We classify the moduli space of vacua by which
terms vanish, and consider branches of vacua. Typical examples are the
Higgs branch $\cM_H$ and the classical Coulomb branch.
The Higgs branch is the symplectic reduction $\bM\tslash G$ explained
in \secref{sec:CH}. It coincides with the hyper-K\"ahler quotient of
$\bM$ by $G_c$ in differential geometry. Connections are trivial ones,
and sections are constant, hence only $\bM$ remain. We do not recall
the notion of hyper-K\"ahler quotients here, so please refer
\cite{MR1193019} for example. Quiver varieties studied by the author
for many years, as well as, toric hyperK\"ahler manifolds are examples
of symplectic reductions, hence of Higgs branches.

On the other hand, the classical Coulomb branch is
$(T^\vee_c\times(\RR^3\otimes\ft_c))/W$, where $T^\vee_c$ is the dual
of a maximal torus $T_c$ of $G_c$, $\ft_c$ is the Lie algebra of
$T_c$, and $W$ is the Weyl group. It is the same as $T^*T^\vee/W$
appeared in \secref{sec:CH}. Sections vanish in the classical Coulomb
branch, and the factor $(\RR^3\otimes\ft_c)$ comes from fields for
which we omit the explanation here. The factor $T^\vee_c$ came from
connections, but they take values in the dual torus $T_c^\vee$ and are
scalars after Fourier transform in an infinite dimensional space of
connections. Even this part of the physics argument is difficult to make
mathematically rigorous, but we will see how $T^\vee$ appears in
\subsecref{subsec:torus} and \thmref{thm:classical} starting from a
mathematically rigorous definition.

Classical Coulomb branches and Higgs branches, and other branches of
the classical moduli space of vacua contain important information of
the supersymmetric gauge theory. It is an important initial step to
analyze the gauge theory. One of the goal of physicists' analysis is a
description of the gauge theory as another supersymmetric quantum
field theory, called a \emph{low energy effective theory}, consisting
of maps from $\RR^3$ with a hyper-K\"ahler manifold as the target
space, together with additional fields, which we will ignore. This is
a surprising statement as field contexts are completely different in
two theories, connections and sections in the original theory while
maps in the new theory. Nevertheless many quantities which physicists
want to compute are the same in two theories in low energy.
The latter super quantum field theory will give Rozansky-Witten invariants after the so-called topological twist.

Hence it is important to determine the target space of the low energy
effective theory. It is roughly the classical moduli space of vacua,
but it is too much hope to expect that local minimum of the lagrangian
contain enough `quantum' information as required in the low energy
effective theory.
Physicists say that the classical Coulomb branch receives quantum
correction. Namely the Coulomb branch $\cM_C$ is
$(T_c^\vee\times(\RR^3\otimes\ft_c))/W$ only in the classical
description, and the actual one is different. It is still a
hyper-K\"ahler manifold as the supersymmetry must exists also in the
low energy effective theory.
This part is difficult to justify directly in mathematically rigorous
way, at least to me. But under our mathematical definition, the
Coulomb branch $\cM_C$ is birational to
$(T_c^\vee\times (\RR^3\otimes\ft_c))/W$, hence $\cM_C$ is a
\emph{correction} of the latter indeed.

Thus the physicists' definition of $\cM_C$ is very far from rigorous
mathematically unlike $\cM_H$. I heard the explanation of the Coulomb
branch in Witten's series of lectures at Newton Institute in 1996
November for the first time, but did not make it as a research object
for many years. Examples of Coulomb branches are familiar
hyper-K\"ahler manifolds to me, hence I had kept interests.

A new insight came to me when I heard Hanany's talk in Warwick in 2014
fall. Hanany explained us there is a formula (monopole formula)
computing the character of the coordinate ring $\CC[\cM_C]$ with
respect to the $\CC^\times$-action. The monopole formula is a sum over
dominant coweights of $G$, and each term is a combinatorial explicit
expression in a coweight. The formula passed many test checking it reproduces the character in many known examples of Coulomb branches.

After Hanany's talk, I looked after a ring whose character reproduces
the monopole formula, as we can reproduce $\cM_C$, at least as an
affine algebraic variety, as the spectrum of the ring. 
Then I found a proposal in \cite{2015arXiv150303676N}, which was
modified in \cite{main}. My path of thinking was explained in
\cite{2015arXiv150303676N}. Let us recall it briefly. The starting
point was \cite[1(iv),(v)]{2015arXiv150303676N}: a hypothetical
topological quantum field theory obtained by a topological twist of
the gauge theory produces a ring: Consider a quantum Hilbert space
$\mathcal H_{S^2}$ for $S^2$. We have a vector in
$\Hom(\mathcal H_{S^2}\otimes\mathcal H_{S^2},\mathcal H_{S^2})$
associated with $M^3$, the $3$-ball with two smaller balls removed
from the interior, which produces a commutative multiplication. Then
the quantum Hilbert space in question is the homology of the moduli
space of solutions of the associated nonlinear PDE on $S^2$, and the
vector is given again by the moduli space of solutions, but on $M^3$
this time, whose image under the boundary value gives a homology
class.
This is an old idea which motivated Atiyah \cite{MR1001453} to write
down axioms of topological quantum field theories based on earlier
works by Donaldson, Floer, and others.

I arrived at a puzzle immediately, as there is only trivial solution
for the nonlinear PDE when $(G_c,\bM) = (\SU(2),0)$, as the only flat
connection on $S^2$ is the trivial one. Since the stabilizer is
nontrivial, namely $\SU(2)$, we may consider the equivariant
cohomology $H^*_{\SU(2)}(\mathrm{pt})$ of a point, but its spectrum is
just $\CC/\pm 1$. It is different from the known answer in physics
(i.e., Atiyah-Hitchin manifold).

I needed a correction, as a naive guess gives an immediate
contradiction.
I made two modifications, (a) forgetting one component of the
nonlinear PDE above, corresponding to the stability condition via the
Hitchin-Kobayashi correspondence, and (b) consider the sheaf of a
vanishing cycle on the moduli space. The latter was motivated by
recent advances in Donaldson-Thomas invariants. 
It will be explained in \secref{sec:not}.
In the joint work \cite{main} I switched from a moduli space on $S^2$
to one on a raviolo\footnote{singular form of \emph{ravioli}, which
  are Italian dumpling.} $\tilde D = D \cup_{D^*} D$, gluing of two
copies of the formal disk $D$ along the punctual disk $D^*$. The
reason was explained in \cite[1(i)]{main}.

\section{A mathematical definition}\label{sec:def}

We will review the definition in \cite{main} in this section.

We assume $\bN$ is a finite dimensional complex representation of a
complex reductive group $G$. Here $\bN$ may not be
irreducible, nor it could be $0$. The symplectic representation $\bM$
is given as $\bN\oplus\bN^*$, but $\bM$ does not show up in this section.

Let $D = \Spec\CC[[z]]$ be the formal disk, $D^\times = \Spec\CC((z))$
the formal punctured disk. We denote $\bN((z))$, $\bN[[z]]$ by
$\bN_\cK$, $\bN_\cO$ respectively. Similarly let $G_\cK = G((z))$,
$G_\cO = G[[z]]$.

The affine Grassmannian $\Gr_G$ is the moduli space
\begin{equation*}
  \left.\left\{ (\scP,\varphi) \middle|
    \begin{aligned}[m]
      & \text{$\scP$ is an algebraic $G$-principal bundle over $D$}\\
      & \text{$\varphi\colon\scP|_{D^\times}\to G\times D^\times$ is a
        trivialization of $\scP$ over $D^\times$}
    \end{aligned}
   \right\}\middle/\text{isom.}\right.
\end{equation*}
It is known that $\Gr_G$ has a structure of an ind-scheme as a direct
limit of projective varieties. Set-theoretically, it is
$\Gr_G = G_\cK/G_\cO$. Namely we take a trivialization of $\scP$ over
$D$ to regard $\varphi$ as an element of $G_\cK$, and kill the
ambiguity of the choice of trivialization by taking the quotient by
$G_\cO$. If we further take the quotient by the left $G_\cO$-action
changing the trivialization $\varphi$, we get
$G_\cO\backslash G_\cK/G_\cO$. It is the moduli space of $G$-bundles
over the raviolo $\tilde D$\footnote{Braverman, my collaborator,
  emphasizes me an importance of use of the raviolo.}.

We then add an algebraic section $s$ of the vector bundle
$\scP\times_G\bN$ associated with the representation $\bN$ to consider
the moduli space $\cT$ of triples
$(\scP,\varphi,s)$. Set-theoretically, it is
$G_\cK\times_{G_\cO}\bN_\cO$. Considering the Taylor expansion of $s$,
we see that $\cT$ is a direct limit of an inverse limit of finite rank
vector bundles over projective schemes. We will consider homology
groups of $\cT$ or its closed varieties, which are rigorously defined
as limit of homology groups of finite dimensional varieties. See
\cite[\S2]{main} for detail.

We introduce a closed subvariety $\cR$ of $\cT$ by imposing the condition that $\varphi(s)$ extends over $D$:
\begin{equation*}
  \cR = \{ (\scP,\varphi,s) \mid \varphi(s)\in\bN_\cO \}/\text{isom.}
\end{equation*}
Since $\varphi$ is a trivialization over $D^\times$, $\varphi(s)$ is
in general has a rational section which may have singularities at the
origin. The space $\cR$ is defined by requiring that coefficients of
negative powers of $\varphi(s)$ are vanishing.
The quotient $G_\cO\backslash\cR$ is the moduli space of pairs of
$G$-bundles and their sections over $\tilde D$.

This space $\cR$ is the main player of our construction. Its meaning
is clearer if we consider a bigger space
\begin{equation*}
  \{ (\scP_1, \varphi_1, s_1, \scP_2, \varphi_2, s_2)
  \in \cT\times\cT \mid \varphi_1(s_1) = \varphi_2(s_2) \}/\text{isom.}
\end{equation*}
This consists of a pair of $G$-bundles over $D$, trivialization over
$D^\times$ and sections of associated vector bundles such that
sections are equal through trivializations. It is a fiber product
$\cT\times_{\bN_\cK}\cT$. If we further require that
$(\scP_2,\varphi_2)$ is the identity element of $\Gr_G$, i.e., the
point where $\varphi_2$ extends across $0\in D$, we recover
$\cR$. Conversely we use the action of $G_\cO$ on $\cR$ to get
$\cT\times_{\bN_\cK}\cT = G_\cK\times_{G_\cO}\cR$ from $\cR$.

From the gauge theoretic view point, $\cT\times_{\bN_\cK}\cT$
parametrizes configurations of a connection and a section on $D$
\emph{twisted} at the origin $0$. Namely $(\scP_1,\varphi_1)$ is
\emph{before} the twist, while $(\scP_2,\varphi_2)$ is
\emph{after}. Since the twisting happens only at the origin, they are
isomorphic outside the origin. Originally we consider a connection and
a section with a point singularity in $2+1$ dimensional space-time in
the $3$-dimensional gauge theory, but we take a $2$-dimensional view
point by looking at two time slices, just before and after the event.

Now the preparation of the space $\cR$ is over, so we consider its
$G_\cO$-equivariant Borel-Moore homology group $H^{G_\cO}_*(\cR)$. We
define its degree so that the fundamental class of the fiber of $\cT$
over the identity element of $\Gr_G$ has degree $0$. We refer
\cite{main} for the precise definition, and omit it here. One can show
that $H^{G_\cO}_*(\cR)$ vanishes in odd degree, and is free over
$H^*_G(\mathrm{pt})$ by using Schubert cell decomposition of the
affine Grassmannian $\Gr_G$.

Next we introduce a convolution product
\begin{equation*}
  \ast\colon H^{G_\cO}_*(\cR)\otimes H^{G_\cO}_*(\cR)
  \to H^{G_\cO}_*(\cR).
\end{equation*}
The rigorous definition in \cite{main} is too technical to be
reproduced here. Let us give a heuristic argument: We formally assume
that we have an induction isomorphism
$H^{G_\cK}_*(\cT\times_{\bN_\cK}\cT)\cong H^{G_\cO}_*(\cR)$, and $\cT$
is smooth. Then using projection to the $(i,j)$-factor
\begin{equation*}
  \cT\times_{\bN_\cK}\cT\times_{\bN_\cK}\cT \xrightarrow{p_{ij}}
  \cT\times_{\bN_\cK}\cT \qquad (i,j) = (1,2), (2,3), (1,3),
\end{equation*}
we define
\begin{equation*}
  c\ast c' = p_{13*}(p_{12}^*c \cap p_{23}^*c').
\end{equation*}
This is not rigorous as we do not know how to define $H^{G_\cK}_*(\cT\times_{\bN_\cK}\cT)$, and $\cT$ is not nonsingular. But we have an alternative rigorous definition of the convolution product $\ast$ on $H^{G_\cO}_*(\cR)$.

We have
\begin{Theorem}
  $(H^{G_\cO}_*(\cR),\ast)$ is a commutative ring.
\end{Theorem}

The method to construct an algebra by convolution has been used in
geometric representation theory, e.g., the group ring of the Weyl
group from Steinberg variety, the universal enveloping algebra of a
Kac-Moody Lie algebra from analog of Steinberg variety for quiver
varieties, etc. But those examples give noncommutative algebras. From
a general theory of convolution, we do not have a reason why $\ast$
becomes commutative.

An explanation of commutativity is given by recalling geometric Satake
correspondence: We consider the abelian category of
$G_\cO$-equivariant perverse sheaves on $\Gr_G$, endow it with a
tensor product via convolution product, and show that the resulted
tensor category is equivalent to one of finite dimensional
representations of Langlands dual group $G^\vee$ of $G$. The latter
category is commutative, i.e., $V\otimes W\cong W\otimes V$, hence the
former is also. A geometric explanation of this commutativity of the
former is given by Beilinson-Drinfeld one-parameter deformation of the
affine Grassmannian. We can give a proof of commutativity in the above
theorem, using this idea \cite{affine}. (In \cite{main} we give
another proof given by a reduction to an abelian case, where it can be
shown by a direct computation.)

Now we get a commutative ring $(H^{G_\cO}_*(\cR),\ast)$. Hence we can
define the affine scheme as its spectrum:
\begin{equation*}
   \cM_C = \Spec (H^{G_\cO}_*(\cR),\ast).
\end{equation*}
We further show that $(H^{G_\cO}_*(\cR),\ast)$ is finitely generated
and integral. Hence $\cM_C$ is an irreducible affine variety. We also
show that it is normal.

A noncommutative deformation is defined as follows. We have a
$\CC^\times$-action on the formal disk $D$ by the loop rotation
$z\mapsto tz$. We have induced actions on various spaces considered
above. In particular, we consider the semi-direct product
$G_\cO\rtimes\CC^\times$ which acts on $\cR$. Hence we can consider
the equivariant Borel-Moore homology group
$H^{G_\cO\rtimes\CC^\times}_*(\cR)$ with respect to the larger group
$G_\cO\rtimes\CC^\times$, and define the convolution product as
above. We thus define \emph{the quantized Coulomb branch} by
\begin{equation*}
   \cAh = (H_*^{G_\cO\rtimes\CC^\times}(\cR),\ast).
\end{equation*}

Convolution products on affine Grassmannian and related spaces were
considered earlier in \cite{MR3013034,MR2135527,MR2422266}, which we
use models for our definition. In \cite{MR3013034}, affine flag
varieties instead of affine Grassmannian, equivariant K-theory instead
of equivariant Borel-Moore homology group were used, but it is
basically understood as a special case of the Coulomb branch where
$\bN$ is the adjoint representation. The algebra constructed there is
Cherednik double affine Hecke algebra (DAHA). If we use affine
Grassmannian instead of flag, we get the spherical part of DAHA. We
get the trigonometric version instead of the elliptic one if we use
homology instead of K-theory. Our Coulomb branch for $\bN=\mathfrak g$
is $\ft\times T^\vee/W$. It is a remarkable example, as the Coulomb
branch does not receive quantum corrections.

In \cite{MR2135527,MR2422266}, the case $\bN=0$ was considered. The
Coulomb branch is the phase space of the Toda lattice for the
Langlands dual group of $G$, or the moduli space of solutions of
Nahm's equation on the interval. We omit further explanation.

\section{Not necessarily cotangent type}\label{sec:not}

In \cite{2015arXiv150303676N} we first made a proposal for the case
when $\bM$ is not necessarily of cotangent type. It was just a
heuristic definition of the coordinate ring $\CC[\cM_C]$ as a graded
vector space, and a definition of the convolution product $\ast$ was
not proposed. Nevertheless another heuristic argument yielded an idea
to define $\CC[\cM_C]$ as $H^{G_\cO}_*(\cR)$ (more precisely homology
of the moduli space on $S^2$). We only have a slight advance in this
direction since \cite{2015arXiv150303676N} was written. Nevertheless
we believe that the original intuition is important, hence we review
it in this section.

The reader can safely skip this section to read other sections.

\subsection{Holomorphic Chern-Simons functional}

Let $\Sigma$ be a compact Riemann surface.
\begin{NB}
We consider a gauge theory
associated with a quaternionic representation $\bM$ of $G$.    
\end{NB}%
We choose and fix a spin structure, i.e., the square root $K_\Sigma^{1/2}$ of
the canonical bundle $K_\Sigma$. We also fix a ($C^\infty$) principal
$G$-bundle $P$ with a fixed reference partial connection
$\DB$.
A {\it field\/} consists of a pair
\begin{itemize}
      \item[] $\DB + A$ : a partial connection on $P$. So $A$ is a
    $C^\infty$-section of $\Lambda^{0,1}\otimes
    (P\times_{G}\mathfrak g)$.
      \item[] $\Phi$ : a $C^\infty$-section of $K_\Sigma^{1/2}\otimes
    (P\times_{G}\bM)$.
\end{itemize}
Let $\mathcal F$ be the space of all fields.
There is a gauge symmetry, i.e., the complex gauge group
$\mathcal G(P)$ of all (complex) gauge transformations of $P$ natural
acts on the space $\mathcal F$.

In fact, as we will see examples below, we need to consider all
topological types of $P$ (classified by $\pi_1(G)$) simultaneously,
but we will ignore this point.

We define an analog of the holomorphic Chern-Simons functional by
\begin{equation}\label{eq:1}
    \CS(A,\Phi)
    = \frac12 \int_\Sigma \omega((\DB + A)\Phi\wedge\Phi),
\end{equation}
where $\omega(\ \wedge\ )$ is the tensor product of the exterior
product and the symplectic form $\omega$ on $\bM$. Since
$(\DB + A)\Phi$ is a $C^\infty$-section of $\Wedge^{0,1}\otimes
K_\Sigma^{1/2}\otimes(P\times_{G}\bM)$, $\omega((\DB +
A)\Phi\wedge\Phi)$ is a $C^\infty$-section of $\Wedge^{0,1}\otimes K_\Sigma
= \Wedge^{1,1}$. Its integral is well-defined.
This is invariant under the gauge symmetry $\mathcal G(P)$.

When $\bM$ is a cotangent type, i.e., $\bM = \bN\oplus\bN^*$, we can
slightly generalize the construction. Let us choose $M_1$, $M_2$ be
two line bundles over $\Sigma$ such that $M_1\otimes M_2 =
K_\Sigma$. We modify $\Phi$ as
\begin{itemize}
      \item[] $\Phi_1$, $\Phi_2$ : $C^\infty$-sections of $M_1\otimes
    (P\times_{G} \bN)$ and $M_2\otimes (P\times_{G}\bN^*)$
    respectively.
\end{itemize}
\begin{NB}
Then $\bmu_\CC(\Phi_1,\Phi_2)$ is a section of $M_1\otimes M_2\otimes (P\times_{G}\mathfrak g^*) = \Lambda^{1,0}\otimes
(P\times_{G}\mathfrak g^*)$.
\end{NB}%
Then
\begin{equation}\label{eq:2}
    \CS(A,\Phi_1,\Phi_2) 
   = \int_\Sigma \langle (\DB+A)\Phi_1,\Phi_2\rangle.
\end{equation}
\begin{NB}
    This is the original version:
\begin{equation*}
    \begin{split}
    & \CS(A,\Phi_1,\Phi_2) 
    = \int_\Sigma \frac12 \left(\langle \DB\Phi_1,\Phi_2\rangle 
    - \langle \DB\Phi_2,\Phi_1\rangle\right)
    + \langle A\wedge\bmu_\CC(\Phi_1,\Phi_2)\rangle
\\
   =\; & \int_\Sigma \langle (\DB+A)\Phi_1,\Phi_2\rangle.
    \end{split}
\end{equation*}
\end{NB}%
It is a complex valued function on $\mathcal F$.

Note that $\mathcal F$ is a complex manifold, in fact, is a complex
affine space, though it is infinite dimensional. Our holomorphic
Chern-Simons functional $\CS$ is a holomorphic function on
$\mathcal F$.

It is easy to see that $(A,\Phi)$ is a critical point of $\CS$ if and
only if the following two equations are satisfied:
\begin{equation}\label{eq:32}
    \begin{split}
        & (\DB+A)\Phi = 0,\\
        & \mu(\Phi) = 0.
    \end{split}
\end{equation}
The first equation means that $\Phi$ is a holomorphic section of
$K_\Sigma^{1/2}\otimes (P\times_{G}\bM)$ when we regard $P$ as a
holomorphic principal bundle by $\DB+A$. The second means that $\Phi$
takes values in $\mu^{-1}(0)$. Therefore $\Phi$ is a holomorphic section
of $K_\Sigma^{1/2}\otimes (P\times_{G}\mu^{-1}(0))$, i.e., a twisted map from $\Sigma$ to the quotient stack $\mu^{-1}(0)/ G$.

Let us denote $\operatorname{crit}(\CS)$ the critical locus of our
holomorphic Chern-Simons functional. Since it is the critical locus of
a holomorphic function on a complex manifold, we could have a sheaf
$\varphi_\CS(\CC_{\mathcal F})$ of vanishing cycle associated with
$\CS$. This is heuristic at this stage as $\mathcal F$ is an infinite
dimensional complex manifold, and hence it is not clear whether the
usual definition of the vanishing cycle can be applied. Nevertheless
it was hoped \cite{2015arXiv150303676N} that one can use an approach
for usual complex Chern-Simons functional for connections on a compact
Calabi-Yau $3$-fold, developed by Joyce and his collaborators
\cite{2012arXiv1211.3259B,2013arXiv1312.0090B}. We thus formally
define
\begin{equation}\label{eq:3}
    H^*_{c,\mathcal G(P)}(\operatorname{crit}(\CS), 
    \varphi_{\CS}(\CC_{\mathcal F}))
\end{equation}
the equivariant cohomology with compact support with the sheaf of
vanishing cycle $\varphi_{\CS}(\CC_{\mathcal F})$ coefficient. The
proposal in \cite{2015arXiv150303676N} was that the dual of this space
(for $\Sigma = \CP^1 = S^2$) has a commutative product, and define the
Coulomb branch as its spectrum.

\subsection{Derived symplectic geometry}

There is an alternative approach for a construction of the perverse
sheaf $\varphi_\CS(\CC_{\mathcal F})$ based on derived symplectic
geometry \cite{MR3090262}, which I learned from Dominic Joyce during a
workshop at Oxford in 2015 after \cite{2015arXiv150303676N} was
written. It is an immediate consequence of results in
\cite{2017arXiv170308578G}. Let us review it for a sake of readers.

Let us first consider $\mu^{-1}(0)/G$ as a derived Artin stack as a
derived fiber product $(\bM/G) \times_{\mathfrak g^*/G} (0/G)$, where
$G$ acts on $\mathfrak g^*$ by the coadjoint action, and
$\bM/G\to \mathfrak g^*/G$ is the moment map. This is equipped with a
$0$-shifted symplectic structure. One of main results in
\cite{MR3090262} is the space $\operatorname{Map}(X,\mu^{-1}(0)/G)$ of
maps from a $d$-dimensional smooth and proper Calabi-Yau $X$ to
$\mu^{-1}(0)/G$ has a $(-d)$-shifted symplectic structure. In
particular, for $\Sigma$ an elliptic curve, (the derived version of)
$\operatorname{crit}(\CS)/\mathcal G(P)$ has a $(-1)$-shifted
symplectic structure when $K_\Sigma^{1/2} = \shfO_\Sigma$.

A modified construction for the case of twisted maps is given in
\cite{2017arXiv170308578G}. It is applicable for our situation of a
compact Riemann surface $\Sigma$. Therefore (the derived version of)
$\operatorname{crit}(\CS)/\mathcal G(P)$ has a $(-1)$-shifted
symplectic structure.

There is an alternative way to define a $(-1)$-shifted symplectic
structure, again due to \cite{2017arXiv170308578G}. We consider the
stack of pairs of $\DB+A$ and $\Phi$ as in \eqref{eq:32}, but without
the equation $\mu(\Phi) = 0$. Let us denote it by
$\operatorname{Sect}_\Sigma(\bM_{K_\Sigma^{1/2}}/G)$. Then the moment
map gives a map to the stack of pairs of $\DB+A$ and $\xi$, a
holomorphic section of $K_\Sigma\otimes (P\times_G \mathfrak g^*)$.
The latter is nothing but the (derived) moduli stack
$\operatorname{Higgs}_G(\Sigma)$ of Higgs bundles, and has a
$0$-shifted symplectic structure. One of main results in
\cite{2017arXiv170308578G} says that the map
\[
\operatorname{Sect}_\Sigma(\bM_{K_\Sigma^{1/2}}/G)\to
\operatorname{Higgs}_G(\Sigma)
\]
is lagrangian. This result was originally observed by Gaiotto
\cite{2016arXiv160909030G} by a heuristic argument as in the previous
subsection. 

There is another lagrangian in $\operatorname{Higgs}_G(\Sigma)$, the
moduli stack $\operatorname{Bun}_G(\Sigma)$ of $G$-bundles on
$\Sigma$. Therefore $\operatorname{crit}(\CS)/\mathcal G(P)$ is a
(derived) fiber product of two lagrangians in a $0$-shifted symplectic
stack, hence has a $(-1)$-shifted symplectic structure by
\cite{MR3090262}.

Now by \cite{2013arXiv1312.0090B} the underlying Artin stack
$\operatorname{crit}(\CS)/\mathcal G(P)$, if it is \emph{oriented},
has a well-defined sheaf of the vanishing cycle, which is regarded as
a definition of $\varphi_\CS(\CC_{\mathcal F})$. We do not recall the
definition of an orientation here, but it is expected that its
existence is guaranteed by the above condition that $\pi_4(G)\to\pi_4(\grpSp(\bM))$ vanishes.

\subsection{Cutting}

Suppose $\bM=\bN\oplus\bN^*$. Then we have a $\CC^\times$-action on
$\mathcal F$ defined by
$t\cdot (A,\Phi_1,\Phi_2) = (A,\Phi_1,t\Phi_2)$. Since $\CS$ is linear
in $\Phi_2$, we have
$\CS(t\cdot(A,\Phi_1,\Phi_2)) = t\CS(A,\Phi_1,\Phi_2)$. Under this
condition for finite dimensional spaces, the vanishing cycle functor
was studied in \cite{2013arXiv1311.7172D}. We hope that this result
can be applied in our infinite dimensional setting, then \eqref{eq:3}
is isomorphic to 
\begin{equation*}
    H^*_{c,\mathcal G(P)}(\cR_\Sigma, \CC),
\end{equation*}
where $\cR_\Sigma$ is the space of $(A,\Phi_1)$ such that
$(\DB+A)\Phi_1 = 0$, that is the space of holomorphic principal
bundles $(P,\DB+A)$ and a holomorphic section of
$M_1\otimes( P\times_G\bN)$. Our space $\cR$ in \secref{sec:def} is
related to $\cR_\Sigma$ by
$G_\cO\backslash \cR = \cR_{\tilde D}/\mathcal G(P)$ though it is not
clear whether we can take $\Sigma = \tilde D$.

\section{Examples}

In order to illustrate that the construction in \secref{sec:def} is
\emph{not} so strange, though we use homology groups of infinite
dimensional spaces, let us give simple examples. This is based on
\cite[\S4]{main}.

\subsection{}\label{subsec:torus}

Let $G = \CC^\times$, $\bN = 0$. This is the simplest case. Since
$\bN=0$, $\cR$ is nothing but the affine Grassmannian $\Gr_G$, and
$\Gr_G$ parametrizes pairs of line bundles on $D$ and their
trivializations over $D^\times$. It is known that $\Gr_G$ with the
reduced scheme structure is the discrete set parametrized by integers
$\ZZ$. In fact, $\varphi(z) = z^n$ is a point corresponding to
$n\in\ZZ$. Therefore
\begin{equation*}
   H^{G_\cO}_*(\cR) = \bigoplus_n H^{\CC^\times}_*(\mathrm{pt}).
\end{equation*}
Note $H^{\CC^\times}_*(\mathrm{pt})$ is the polynomial ring $\CC[w]$
of one variable $w$. Since we have a polynomial ring over each integer
$n$, we need to calculate the product of a polynomial on $m$ and that
on $m$. Since we do not give the precise definition of the convolution
product, we cannot perform the check, but for $G=\CC^\times$, the
product $\ast$ is given by the push-forward homomorphism of the map
given by tensor product
\begin{equation*}
   \Gr_{\CC^\times}\times\Gr_{\CC^\times}\xrightarrow{\otimes}\Gr_{\CC^\times}.
\end{equation*}
Then the product of $f(w)$ on $m$ and $g(w)$ on $n$ is $f(w)g(w)$ on
$m+n$. Let us denote by $x$ the polynomial $1$ on the integer
$n=1$. We then have
\begin{equation*}
   H^{G_\cO}_*(\cR) \cong \CC[w,x^\pm] = \CC[\CC\times\CC^\times].
\end{equation*}
Therefore the Coulomb branch is $\CC\times\CC^\times$. Since this is
nothing but $\RR^3\times S^1$, the Coulomb branch does not receive the
quantum correction. This is a reflection of the fact that the gauge
theory is trivial in this case.

Let us further consider the case when $G$ is a torus $G$, and
$\bN = 0$. Then $\Gr_T$ is a discrete space parametrized by
$\Hom(\CC^\times,T)$. Therefore
$H^{T_\cO}_*(\cR) = \bigoplus_{\lambda\in\Hom(\CC^\times,T)}
H^*_T(\mathrm{pt})$. Note that $H^*_T(\mathrm{pt})$ is the space
$\CC[\ft]$ of polynomials on the Lie algebra $\ft$ of $T$.  On the
other hand, let $e^\lambda$ denote the fundamental class of the point
$\lambda$. We have $e^\lambda\ast e^\mu = e^{\lambda+\mu}$ as above.
Since this can be regarded as the ring of characters of the dual
$T^\vee$ of $T$, the Coulomb branch is
$\ft\times T^\vee = T^* T^\vee$.

\subsection{}\label{subsec:C2}

Let us keep $G$ as $\CC^\times$, and replace the representation to the
standard on $\bN = \CC$.
As we have already explained, $\Gr_{\CC^\times}$ is a discrete set
parametrized by $\ZZ$, and $\cR$ consists of vector spaces over
integers $n\in\ZZ$. Since the condition is that we do not get
singularities by $\varphi(z) = z^n$, we have
\begin{equation*}
  \cR = \bigsqcup_{n\in\ZZ} z^{n}\CC[z]\cap \CC[z] =
  \bigsqcup_{n\in\ZZ} z^{\max(0,n)}\CC[z].
\end{equation*}
By the Thom isomorphism for each $n$, we have
$H^{G_\cO}(\cR) \cong \bigoplus_n
H^{\CC^\times}_*(\mathrm{pt})$. Hence it is the same as above example
as a vector space. On the other hand, the convolution product is
different. In fact, products of homology classes over $n>0$ and those
over $n < 0$ are different from above. We cannot check the assertion
as we omit the definition, but the product of the fundamental classes
of $n=1$ and $n=-1$ is the image under the pushforward homomorphism for 
\begin{equation*}
  z\CC[z] \to \CC[z]
\end{equation*}
of the fundamental class. Since the image of this map is a codimension
$1$ subspace, it is nothing but the cup product of $w$ with the
fundamental class. Therefore if we denote fundamental classes of
$n=1$, $-1$ by $x$, $y$ respectively, we get $xy = w$. Thus
\begin{equation*}
  H^{G_\cO}_*(\cR) \cong \CC[w,x,y]/(w=xy) \cong \CC[x,y] = \CC[\CC^2].
\end{equation*}
Namely the Coulomb branch in this case is $\CC^2$.

If we replace the representation by the $1$-dimensional representation
with weight $N$, the product $xy$ is replaced by the image of the
fundamental class under $z^{|N|}\CC[z]\to \CC[z]$. Therefore the
coordinate ring is $\CC[w,x,y]/(w^{|N|}=xy)$. Hence the Coulomb branch
is nothing but the simple singularity of type $A_{|N|-1}$.

\section{Structures}

In this section we review several structures of the Coulomb branch
$\cM_C$. We also discuss the corresponding structures for the Higgs
branch $\cM_H$. They have been discussed in physics context. A point
is that they can be realized rigorously in the definition in
\secref{sec:def}.

\subsection{} (See \cite[\S4(iii)(a)]{2015arXiv150303676N} and
\cite[Remark 2.8]{main}.)
$H^{G_\cO}_*(\cR)$ is a graded algebra by the half of homological
degree. We thus have a decomposition
$\CC[\cM_C] = \bigoplus_d \CC[\cM_C]_d$ such that
$\CC[\cM_C]_d\cdot\CC[\cM_C]_{d'}\subset\CC[\cM_C]_{d+d'}$. It means
that $\cM_C$ has a $\CC^\times$-action. In fact, $\CC[\cM_C]_d$ is the
weight space with respect to the $\CC^\times$-action with weight $d$.

In above examples, the $\CC^\times$-actions are weight $1$ on $x$, and
$0$ on $y$. Thus they are the standard action on the first factor of
$\CC\times\CC^\times$ and $\CC^2 = \CC\times\CC$ respectively.

Remark that in general, degrees take values in integers, not
necessarily nonnegative. Therefore $\cM_C$ may not be cone. Here
$\cM_C$ is a cone if $\CC[\cM_C]_d = 0$ ($d < 0$), $\CC[\cM_C]_0 = \CC$.

In physics context, it is expected that the $\CC^\times$-action, or
rather the $S^1$-action, extends to an $\SU(2)$-action after a certain
correction. We do not explain the correction, but it is given by a
hamiltonian torus action explained below. In particular, there will be
no correction when $G$ is semisimple. The induced $\SU(2)$-action on
the two sphere of complex structures
$S^2 = \{ a I + bJ + cK \mid a^2 + b^2 + c^2 = 1\}$ is the standard
one through $\SU(2)\to\SO(3)$, where $(I, J, K)$ is the hyper-K\"ahler
structure. Once we fix a complex $I$, we could see only the
$S^1$-action preserving $I$. This is the reason why we could only see
the $S^1$-action in the current definition, which does not realize the
hyper-K\"ahler structure.

For example, we have an $\SU(2)$-action on
$\CC\times\CC^\times = \RR^3\times S^1$, once we view $\RR^3$ as
$\su(2)$. Our $S^1$-action has the half weight. For $\CC^2$, we
correct the action by a hamiltonian $S^1$-action with weights $-1/2$ on
$x$, $1/2$ on $y$. If we multiply weights by two, it becomes the
restriction of the standard $\SU(2) = \grpSp(1)$-action, given by the
identification $\CC^2$ with the quaternion field $\mathbb H$. (It is
not a complex linear, hence it is different from the standard
$\SU(2)$-action on $\CC^2$. They are left and right multiplication of quaternions respectively.

Let us consider the Higgs branch $\cM_H$ where the $\SU(2)$-action can
be easily described. The quaternionic vector space $\bM$ has an
$\SU(2)=\grpSp(1)$-action by multiplication of quaternion. It commutes
with the $G$-action, hence we have an $\SU(2)$-action on $\cM_H$. It
rotates the two sphere $S^2$ of complex structures, as it is so on
$\bM$.

\subsection{} (See \cite[\S3(vi)]{main}.)
Since $H^{G_\cO}_*(\cR)$ is an equivariant homology group, there is a
homomorphism from $H^*_{G_\cO}(\mathrm{pt})\cong
H^*_G(\mathrm{pt})$. (Remark that the convolution product $c\ast c'$
is not \emph{naturally} $H^*_G(\mathrm{pt})$-linear, in fact it
\emph{isn't} on the noncommutative deformation.)

Taking spectrum, we obtain
\begin{equation*}
    \varpi\colon\cM_C\to \Spec H^*_G(\mathrm{pt}).
\end{equation*}
It is well-known that
\begin{equation*}
     H^*_G(\mathrm{pt}) = \CC[\mathfrak g]^G
    = \CC[\ft]^W,
\end{equation*}
and hence $\Spec H^*_G(\mathrm{pt}) = \ft/W$, where
$\ft = \operatorname{Lie}T$. This is an affine space.

This construction remains on the noncommutative deformation:
\begin{equation*}
    H^{G\times\CC^\times}_*(\mathrm{pt})\to
    \cAh = H^{G_\cO\rtimes\CC^\times}(\cR).
\end{equation*}
This is an injective algebra homomorphism. In particular, the
noncommutative deformation $\cAh$ contains a large commutative
subalgebra. Considering the specialization at $\hbar=0$, we deduce
that $\varpi$ is Poisson commuting. Namely pull-backs of functions
$f$, $g$ on $\ft/W$ satisfy $\{ \varpi^*f, \varpi^*g\} = 0$.

We have the following
\begin{Theorem}[See \protect{\cite[\S5(v)]{main}}.]\label{thm:classical}
  A generic fiber of $\varpi$ is $T^\vee$. More precisely we have the following commutative diagram, whose upper horizontal arrow is birational:
\begin{equation*}
    \xymatrix{
    \cM_C \ar@{.>}[rr] \ar[dr]_\varpi &&
    T^* T^\vee/W = \ft\times T^\vee/W \ar[dl]^{\ \ \ \text{the first projection}}
    \\
    & \ft/W &}
\end{equation*}
\end{Theorem}

This is a consequence of the fixed-point localization theorem for the
equivariant homology group. The localization theorem says that we have an isomorphism
\begin{equation*}
    H^{T_\cO}_*(\cR)\otimes_{H^*_T(\mathrm{pt})}\mathbb F
    \cong H^{T_\cO}_*(\cR^T)\otimes_{H^*_T(\mathrm{pt})}\mathbb F,
\end{equation*}
where $\mathbb F$ is the quotient field of $H^*_T(\mathrm{pt})$. Here
$\cR^T$ is the set of $T$-fixed points in $\cR$, and the isomorphism
is the pushforward homomorphism of the inclusion
$\cR^T\hookrightarrow\cR$. Combining this with the fact that
$H^{G_\cO}_*(\cR)$ is the $W$-invariant part of $H^{T_\cO}_*(\cR)$, it
becomes enough to compute the equivariant homology group of
$\cR^T$. Since $\cR^T$ is $\Gr_T\times\bN^T$, the calculation in
\subsecref{subsec:torus} shows that it is $\ft\times T^\vee$.

The operation $\otimes_{H^*_T(\mathrm{pt})}\mathbb F$ corresponds to
the restriction of the generic point of $\ft/W$. This is a standard
argument which tells that it is useful to view equivariant homology
groups as families over $\ft/W$.

In conclusion, $\varpi$ is Poisson commuting and has algebraic tori as
fibers. Hence $\varpi\colon\cM_C\to\ft/W$ is an integrable system in
the sense of Liouville, and $\cAh$ is its quantization.

For the Higgs branch $\cM_H$, we do not have a general construction of
an integrable system, though we could see it in many examples.

\subsection{}\label{subsec:hamtori} (See \cite[\S4(iii)(c)]{2015arXiv150303676N} and
\cite[\S3(v)]{main}.)
It is known that the affine Grassmannian $\Gr_G$ is topologically a
based loop group $\Omega G$. In particular, its connected components
are in bijection to the fundamental group $\pi_1(G)$ of $G$. It is
well-known that $\pi_1(G)$ is a finitely generated abelian group. The
homology group of $\cR$ decomposes according to connected components
of $\cR$, which are the same as those of $\Gr_G$. This decomposition
is compatible with the convolution product: let $\cR_\gamma$ denote
the connected component corresponding to $\gamma\in\pi_1(G)$. Then we
have
$H^{G_\cO}_*(\cR_\gamma)\ast H^{G_\cO}_*(\cR_{\gamma'})\subset
H^{G_\cO}_*(\cR_{\gamma+\gamma'})$.

In terms of $\cM_C = \Spec H^{G_\cO}_*(\cR)$, this decomposition means
that the Pontryagin dual $\pi_1(G)^\wedge = \Hom(\pi_1(G),\CC^\times)$
of $\pi_1(G)$ acts on $\cM_C$. In above examples, we have
$\pi_1(G) = \pi_1(\CC^\times) = \ZZ$, and its Pontryagin dual is
$\CC^\times$. The action is on the second factor in the first example
$\cM_C=\CC\times\CC^\times$. In the second example, $x$ has weight $1$
and $y$ has weight $-1$.

Since this action extends to the noncommutative deformation
$H^{G_\cO\rtimes\CC^\times}_*(\cR)$, it follows that the symplectic
form is preserved under the action.

When $G$ is semisimple, $\pi_1(G)$ is a finite group, and its
Pontryagin dual also. We obtain a torus when $\Hom(G,\CC^\times)$ is
nontrivial. Let $\chi\in\Hom(G,\CC^\times)$. The corresponding moment
map of the $\CC^\times$-action via
$\Hom(G,\CC^\times)\cong\Hom(\CC^\times,\pi_1(G)^\wedge)$ is given by
the composition of $\varpi$ with
$d\chi\colon \mathfrak g\to\operatorname{Lie}\CC^\times$. In
particular, the action is hamiltonian. One can also show that the
symplectic reduction of $H^{G_\cO\rtimes\CC^\times}_*(\cR)$ is the
Coulomb branch of the kernel of $\chi$. See \cite[\S3(vii)(d)]{main}.

For the Higgs branch, $\chi\in\Hom(G,\CC^\times)$ is used to introduce
a stability condition for the geometric invariant theory
quotient. Namely we can consider $\operatorname{Proj}$ of the graded ring
$\bigoplus_{n=0}^\infty \CC[\mu^{-1}(0)]^{G,\chi^n}$ of
semi-invariants. Here
$\CC[\mu^{-1}(0)]^{G,\chi^n} = \{ f\in\CC[\mu^{-1}(0)] \mid f(g\cdot
x) = \chi(g)^n f(x)\}$.
Also we can use
$\zeta\in\Hom(\mathfrak g,\operatorname{Lie}\CC^\times)$ to
perturb the defining equation as $\mu = \zeta$.

\subsection{} (See \cite[\S5(i)]{2015arXiv150303676N} and
\cite[\S3(viii)]{main}.)
Suppose that $\bN$ is a representation of a larger group
$\widetilde G$ containing $G$ as a normal subgroup. The quotient group
$\widetilde G/G$ is called the \emph{flavor group} in physics
literature. Let us denote it by $G_F$.

Since $\widetilde G_\cO$ acts on $\cR$, we can consider the
equivariant homology group $H^{\widetilde G_\cO}_*(\cR)$ with respect
to the larger group $\widetilde G_\cO$. It is a commutative ring over
$H^*_{G_F}(\mathrm{pt})$, hence the corresponding spectrum is a family
of varieties over
$\Spec H^*_{G_F}(\mathrm{pt}) = \Spec
\CC[\mathfrak g_F]^{G_F}$.
The fiber over $0$ is the original $\cM_C$. Namely $\cM_C$ has a
deformation parametrized by $\mathfrak g_F\dslash G_F$.

Although we omit the detail, we can construct (candidates) of partial
resolutions of $\cM_C$ corresponding to cocharacters of a maximal
torus $T_F$ of $G_F$. See \cite[\S3(ix)]{main}.

On the Higgs branch $\cM_H$, we have an induced action of
$G_F = \widetilde G/G$. Note that structures in this and previous
subsections are swapped for $\cM_C$ and $\cM_H$. Namely
$\Hom(G,\CC^\times)$ gives a deformation/resolution on $\cM_H$ and a
group action on $\cM_C$. On the other hand $G_F$ gives a group action
on $\cM_H$ and a deformation/resolution on $\cM_C$.

\subsection{}

Let us consider toric hyper-K\"ahler manifolds as examples of
structures of one and two subsections before. We start with an exact
sequence of tori
\begin{equation*}
    1 \to T = (\CC^\times)^{d-n} \to \widetilde T = (\CC^\times)^d
    \to T_F = (\CC^\times)^n \to 1.
\end{equation*}
We take the standard representation $\bN = \CC^d$ of $\widetilde T$
and denote its restriction to $T$ also by $\bN$. We have
$\cM_C(\widetilde T,\bN)\cong \CC^{2d}$ by the computation in
\subsecref{subsec:C2}. By the construction of two subsections before,
the Pontryagin dual of $\pi_1(\widetilde T)$ acts on $\CC^{2d}$. This
is nothing but the standard action of the dual torus
$\widetilde T^\vee$ of $\widetilde T^\vee$. The dual $T_F^\vee$ of
$T_F$ is a subtorus of $\widetilde T^\vee$, hence acts on
$\CC^{2d}$. As we explained in two subsections before, the Coulomb branch $\cM_C(T,\bN)$ for the subgroup $T$ is nothing but the symplectic quotient $\CC^{2d}\tslash T_F^\vee$ of $\CC^{2d}$ by $T_F^\vee$. The exact sequence of dual tori
\begin{equation*}
    1\to T_F^\vee \to \widetilde T^\vee \to T^\vee \to 1
\end{equation*}
identifies it as the Higgs branch for $T_F^\vee$ for the
representation $\CC^d$.  Namely under the exchange
$T\leftrightarrow T_F^\vee$, the Higgs and Coulomb branches are
exchanged.

\subsection{} (See \cite[\S4(iii)(d)]{2015arXiv150303676N} and
\cite[App.~A]{affine})

We can extend the hamiltonian torus action from $\Hom(G,\CC^\times)$
to a nonabelian group action sometimes. Suppose that we have a
subspace $\mathfrak l$ in $\CC[\cM_C]$ which is a Lie subalgebra with
respect to the Poisson bracket $\{\ ,\ \}$. For example, the space of
degree $1$ elements forms a Lie subalgebra as the Poisson bracket is
of degree $-1$. We consider Hamiltonian vector fields $H_f$ for
$f\in \mathfrak l$, and they form a Lie subalgebra in the Lie algebra
of vector fields on $\cM_C$ as $[H_f, H_g] = H_{\{f,g\}}$. Thus
$\mathfrak l$ acts on $\cM_C$ so that the transpose of the moment map
is is the natural homomorphism
$\CC[\mathfrak l^*] = \operatorname{Sym}(\mathfrak l)\to \CC[\cM_C]$.
In many examples, $\mathfrak l$ is integrated to a Lie group action.

Consider the example~\ref{subsec:torus}. The symplectic form, in this
example, is a standard one $dw\wedge\frac{dx}x$. We have
$\{ x, w\} = w$, and $\CC x\oplus \CC w$ is a $2$-dimensional Lie
subalgebra. This is integraded to a $\CC^\times\ltimes \CC$-action as
$(t,s)(x,w) = (tx, sx + w)$ for $(t,s)\in\CC^\times\ltimes \CC$.

This computation is not enlightening as we know the Coulomb branch
explicitly. One can consider also the example~\ref{subsec:C2}, but
again not enlightening.
A nontrivial example is the action of $\operatorname{Stab}_{G_Q}(\mu)$
on a slice to ${\Gr}_{G_Q}^\mu$ in $\overline{\Gr}_{G_Q}^{\la}$ as the
Coulomb branch of a quiver gauge theory explained in the next section.
See \cite[App.~A]{affine}.

\section{Quiver gauge theories}

At the time of this writing, Coulomb branches of $(G,\bN)$ whose Higgs
branches are quiver varieties are most studied. Let $Q$ be a quiver
with the vertex set $Q_0$ and the edge set $Q_1$. For an edge
$h\in Q_1$, let denote the starting and ending vertices by $\vout{h}$,
$\vin{h}$ respectively. For given two $Q_0$-graded finite dimensional
complex vector spaces $V = \bigoplus V_i$, $W = \bigoplus W_i$, we set
\begin{equation*}
    \begin{split}
        & G = \prod_{i\in Q_0}\GL(V_i),
        \\
        & \bN = \bigoplus_{h\in Q_1}\Hom(V_{\vout{h}},V_{\vin{h}})
        \oplus\bigoplus_{i\in Q_0} \Hom(W_i, V_i).
    \end{split}
\end{equation*}
The pair $(G,\bN)$ is a quiver gauge theory. Here the $G$-action on
$\bN$ is the natural one.

When $Q$ is of type $ADE$, the Coulomb branch $\cM_C$ is identified
with a moduli space of monopoles on $\RR^3$ with singularities at the
origin in physics. This assertion is proved in the above mathematical
definition when the monopole moduli space is replaced by its
algebro-geometric analog (\cite{2016arXiv160403625B}). Here the
structure group of monopoles is the complex simple Lie group $G_Q$ of
type $Q$ of the adjoint type, the dimensions of $V_i$ give the charges
of monopoles, and the dimensions of $W_i$ determine the singularity
type.
The definition of algebro-geometric analog is not simple in general,
but when $\mu = \sum \dim W_i \varpi_i - \dim V_i\alpha_i$ is
dominant, it is given as follows: Consider the affine Grassmannian for
$G_Q$, and Schubert varieties $\overline{\Gr}_{G_Q}^{\la}$,
$\overline{\Gr}_{G_Q}^{\mu}$ for $\lambda = \sum \dim W_i\varpi_i$ and
$\mu$. Then the intersection of a transversal slice to
${\Gr}_{G_Q}^\mu$ and $\overline{\Gr}_{G_Q}^{\la}$ is
$\cM_C$.

Under the geometric Satake correspondence, the affine Grassmannian is
connected with representation theory of the Langlands dual group
$G_Q^\vee$ of $G_Q$. On the other hand, homology groups of quiver
varieties have structures of representations of the Lie algebra of
$G_Q$, or of $G_Q^\vee$ which is the simply-connected type.
The symplectic duality mentioned in the introduction is (and should
be) formulated so that two constructions are related by a `duality'.

To determine Coulomb branches, we use the following recipe:
\begin{enumerate}
\item First construct a candidate of $\cM_C$. In many cases, we just
  take an answer given by physicists.
\item Next construct an integrable system on the candidate, which is
  expected to correspond to $\varpi$.

\item Show that the integrable system is a flat family, and $\cM_C$ is
  normal.

\item The birational isomorphism between $\cM_C$ and the candidate
  through $T^*T^\vee/W$ extends over the complement of the inverse
  image of a codimension $2$ subset in $\ft/W$.
\end{enumerate}
It is a consequence of the normality that the extension outside
codimension $2$ guarantees the isomorphism everywhere.
As we explained above, $\cM_C$ is birational to $T^*T^\vee/W$ by an
application of the localization theorem in equivariant homology
groups. By a similar argument, $\cM_C$ can be determined at a
codimension $1$ subvariety by a reduction to Coulomb branches of tori
and rank $1$ groups. The abelian cases are determined as in
\subsecref{subsec:C2}, and the rank $1$ case is a hypersurface in
$\CC^3$ (\cite[\S6(iv)]{main}). Therefore (4) is usually an easy
step. On the other hand (3) is checked by a case-by-case argument, is
usually key point of the proof.

When $Q$ is affine type $ADE$, we replace monopoles by instantons. We
should consider instantons on the Taub-NUT space, not on $\RR^4$ in
general. When $\mu$ is dominant, it is expected that moduli spaces on
$\RR^4$ and on the Taub-NUT space are isomorphic as complex symplectic
varieties. (Hyper-K\"ahler metrics are different.)

For instanton moduli spaces, either on $\RR^4$ or the Taub-NUT space,
the property (3) is not known. Hence we cannot prove that Coulomb
branches are instanton moduli spaces in general.

In fact, (3) is a delicate property. For example, nilpotent orbits are
normal for type $A$, but not in general. On the other hand, Coulomb
branches are always normal. It is known that nilpotent orbits and
their intersection with Slodowy slices for classical groups appear as
Higgs branches. A naive guess gives the corresponding Coulomb branches
are also, but they should not by the normality. Hanany et al find
examples of Coulomb branches, which are normalization of non-normal
nilpotent orbits.

For affine type $A$, we can use Cherkis bow varieties instead of
instanton moduli spaces on the Taub-NUT space. Bow varieties are
moduli spaces of solutions of Nahm's equation, which is a nonlinear
ODE. The ODE is hard to analyze, hence we rewrite bow varieties by
moduli spaces of representations of a quiver with relations, and show
the property (3) (see \cite{2016arXiv160602002N}). Thus Coulomb
branches for affine quiver gauge theories of type $A$ are all
determined.

\section{Quantized Coulomb branches}

Less is known for quantized Coulomb branches than Coulomb branches
themselves.

For a quiver gauge theory of finite type $ADE$, the quantized Coulomb
branch $\cAh$ is isomorphic to a shifted Yangian, as proved in
Appendix of \cite{2016arXiv160403625B}. But this was shown under the
assumption that $\mu$ is dominant. General cases remain open.

We have mentioned that the quantized Coulomb branch for
$\bN = \mathfrak g$ is the spherical DAHA. Consider the case
$G=\GL(k)$ as an example of a quiver gauge theory for the Jordan
quiver with $V=\CC^k$, $W=0$. We generalize this case to $V=\CC^k$,
$W=\CC^r$. In this case $\cAh$ is the spherical part of the rational
Cherednik algebra associated with the wreath product
$\ZZ/r\ZZ\wr S_k = (\ZZ/r\ZZ)^k\rtimes S_k$
\cite{2016arXiv160800875K}. (The corresponding Coulomb branch is
$\operatorname{Sym}^k(\CC^2/(\ZZ/r\ZZ))$.)

\bibliographystyle{myamsalpha}
\bibliography{nakajima,mybib,coulomb}    

\end{document}